%% file: review_semi_class.tex
\documentclass[11pt,twocolumn,a4paper]{article}

\usepackage{latexsym}
\usepackage{amsmath}
\usepackage{amssymb}
\usepackage{amsthm}
\usepackage{graphics}
\usepackage{epsf}
\usepackage{hyperref}
\usepackage{mathrsfs}
\usepackage{xspace}
\usepackage{leftidx}

\input{macros}

\newtheorem*{thm}{\sc Theorem}

\author{Clotilde Fermanian Kammerer \& J\'er\^ome Le Rousseau}

\title{Semi-classical analysis}

\begin{document}

\maketitle


\section{Introduction}

\subsection{What is semi-classical analysis ? }

Semi-classical analysis has its roots in the foundations of quantum
mechanics. Simultaneously with this new theory
arose the question of understanding the links between classical and
quantum mechanics.  It turned out that the Planck constant $\hbar$ can
be understood as the obstruction to give a classical description of a
quantum particule by the simultaneous knowledge of its position and
its momentum. This is expressed by the Heisenberg uncertainty principle
that we first discuss.  \smallskip

In quantum mechanics, a particule is described by a probability
measure $|\psi(x)|^2 d x$, with $\psi$
a normalized square integrable function
on the configuration space $\R^d_x$, called its
wave function. Denoting by $x_j$ the coordinates of $x\in\R^d$, the average position
of the particule is 
\[
\jpsi{x_j}=\int_{\R^d} x_j |\psi(x)|^2 dx,\; 1\leq j\leq d,
\]
that is, 
the expectation
value of the observable $x_j$.
Similarly, the average momentum is 
\begin{equation}\label{def:av_imp}
\jpsi{\xi_j}=\int_{\R^d}  \hbar D_{x_j}\psi(x) \, \overline\psi(x)  dx,\;\; D_{x_j}=\frac 1i\partial_{x_j}.
\end{equation}
Considering  the variance of these random variables,
\begin{align*}
&(d_\psi x_j)^2=  \bjpsi{ (x_j- \jpsi{x_j})^2},\\
&(d_\psi \xi_j)^2= \bjpsi{(\xi_j- \jpsi{\xi_j})^2},
\end{align*}
 the  Heisenberg uncertainty principle  reads
\begin{equation*}
d_\psi x_j\, d_\psi \xi _j\geq \frac \hbar 2,\;\;1\leq j\leq d.
\end{equation*}
It relies on the Cauchy-Schwarz inequality
\begin{align*}
  &\left|\Im \biginp{ (x_j- \jpsi{x_j})\psi}{(\hbar D_{x_j} -
    \jpsi{\xi_j})\psi}_{\! L^2}\right|\\
  &\quad \leq  \| (x_j- \jpsi{ x_j})\psi\|_{L^2} \|( \hbar D_{x_j}
    - \jpsi{\xi_j})\psi\|_{L^2}\\
  &\quad= d_\psi x_j\, d_\psi \xi _j,
\end{align*}
and the observation 
\begin{align*}
&\Im \biginp{ (x_j- \jpsi{x_j})\psi}{(\hbar D_{x_j} -
  \jpsi{\xi_j})\psi}_{\! L^2}\\
  &\quad = \frac 1{2i}
    \biginp{\! \left[\hbar D_{x_j} - \jpsi{\xi_j}, x_j- \jpsi{x_j}\right]
    \psi}{\psi}_{\! L^2}
 =-\frac \hbar 2.
\end{align*}

The Planck constant $\hbar$ reflects the difference between quantum
and classical mechanics, since, in the latter, the position and the
momentum are deterministic variables. The subject of semi-classical
analysis is to understand how one can derive classical mechanics from
quantum mechanics.
Even though $\hbar$ is a physical constant, this
is done by performing the limit $\hbar\rightarrow 0$.  For this
reason, we will skip the notation $\hbar$ and denote by $h$ a small
parameter that is present in some problems of interest involving PDEs.
Carrying a semi-classical analysis of this problem consists in investigating the properties of a phenomenon of interest in the limit $h\rightarrow 0$. 
This type of analysis led to  the 
development of asymptotic technics that are now used in various fields of applied mathematics. 
Examples are the determination of the asymptotics of the spectrum of
Schr\"odinger operators or the characterization of the properties of the
solutions to time-dependent Schr\"odinger equations.

\subsection{Outline}
We introduce in
Section~\ref{sec: Some semi-classical problems} three representative topics in semi-classical analysis.  Starting from the
correspondence between classical and quantum mechanics, basic
semi-classical analysis tools and results are presented in
Section~\ref{sec:cor_pri}. In Section~\ref{sec:applications},
 the three problems of
Section~\ref{sec: Some semi-classical problems} are investigated in
the light of the introduced techniques allowing one to emphasize
different aspects of semi-classical analysis.

\section{Some semi-classical problems}
\label{sec: Some semi-classical problems}

Three problems are presented. They originate from various fields: theoretical
chemistry, spectral geometry, and control theory. 
In each case the semi-classical parameter has a different
interpretation. 

\subsection{Schr\"odinger equation in the Born-Oppenheimer approximation}

The dynamics of a molecule consisting in $k_e$ electrons and $k_n$ nuclei of masses $(M_j)_{1\leq j,\leq k_n}$ (in atomic units)  is described by a 
 wave function belonging to  $L^2(\R^{3k_e+3k_n})$. 
Dating from the 30s, the {\it Born-Oppenheimer approximation}~\cite{BOp}  suggests to take advantage of the fact that, $m_e$ being the mass of an electron, the ratio $m_e/M_j$ is small, for all the nuclei, and roughly, of the main size, even though the $j$-ths atoms are different. Setting 
$$\sqrt{\frac{m_e}{M_j}}\sim h,\;\; 1\leq j\leq k_n,$$
one introduces in the equations the small parameter $h$ and writes 
$$\widehat H_{\rm mol}  =  - \frac{h^2}{2} \Delta_{x}+\widehat H_e(x),$$
where $x$ is in $\R^{3k_n}$ and denotes the coordinates of the nuclei and the electronic Hamiltionian $\widehat H_{\rm e}(x)$ takes into account the kinetics of the electrons, together with the interactions between the electrons themselves, nuclei, and electron/nuclei. 
\smallskip

For all $x$ in $\R^{3k_n}$, the operator $\widehat H_{\rm e}(x)$ is a self-adjoint operator on $L^2(\R^{3k_e})$ with spectrum $\sigma_{\rm e} (x)$ that depends on the configuration $x$ of the nuclei. When the initial data $\psi^h_0$ is in the vector-sum of $N$ eigenspaces  of $H_{\rm e}(x)$ corresponding to $N$ eigenvalues isolated from the remainder of the spectrum, it 
 has been proved in~\cite{ST,MS}, that, considering semi-classical times $t\sim \frac 1h$, one is left with a system of semi-classical Schr\"odinger equations 
\begin{equation}\label{eq:schro_BO}
ih  \partial_t \psi^h = - \frac {h^2}{2}\Delta  \psi^h +V(x)\psi^h,\;\; (t,x)\in \R\times \R^d,
\end{equation}
with $\psi ^h\in L^2(\R^d, \C^{N})$ and $V$ a smooth matrix-valued
potential.  The analysis is thus reduced to a finite number of spectral
components, and, as discussed in Section~\ref{sec:app_BO},
semi-classical technics allow one to develop numerical tools adapted for solving these equations~\cite{LaLu}.

\subsection{Eigenfunctions of the Laplacian and quantum limits}

Let us consider $(M,g)$ a smooth compact Riemannian manifold without boundary. The  Laplace-Beltrami operator  $-\Delta_M$ is a nonnegative self-adjoint operator with compact resolvent, and admits a 
sequence of normalized eigenfunctions  $(\varphi_k)_{k\in\N}$ and eigenvalues $(E_k)_{k\in\N}$, ordered in increasing order: 
\begin{align}\label{eq:eigen}
&- \Delta_M \varphi_k=E_k\varphi_k,\\
\nonumber
&\; 0=E_1\leq E_2\leq  \cdots \leq E_k\leq \cdots,\;\; E_k\Tend{k}{\infty} +\infty.
\end{align}
A historical question~\cite{Einstein} concerns   the densities 
\[
\nu_k(x)=|\varphi_k(x)|^2 dx,
\]
and the analysis of their  limit points, measures on $M$, as $k\to+\infty$.
Such measures are called {\it quantum limits}.
Setting 
\[
h_k=\frac 1{\sqrt E_k},
\]
one is left with a semi-classical problem consisting in the analysis
of a sequence of wave functions $(\varphi_k)_{k\in\N}$ satisfying the
semi-classical PDE
\begin{equation*}
-h_k^2\Delta_M \varphi_k = \varphi_k.
\end{equation*}
As we shall see in Section~\ref{sec:app_geo}, this approach of the problem allows one to derive fundamental
properties of the quantum limits, leading in certain cases, to their
determination (see the {\it Schnirelman Theorem} and its proofs by Y.~Colin de Verdi\`ere and S.~Zelditch, independently,~\cite{Sni,CdV:85,Zel87}, or the surveys~\cite{AFF,Anantharaman_book}). 
\smallskip

This type of question is also posed in the context of random surfaces
with genus that tends to infinity, the semi-classical parameter 
is then the inverse of the genus~\cite{Monk}. These examples and the
preceding one illustrate that the physical  meaning of the semi-classical
parameter may be far from the actual Planck constant~$\hbar$.

\smallskip 
In the preceding two
examples, the small scale $h$ appears naturally and its presence in
the equations endows the solutions with specific features. For example, the family of eigenfunctions  $(\varphi_k)_{k\in\N}$ in~\eqref{eq:eigen} have $H^s$-Sobolev norms of  size $h_k^{-s}$. One can also argue in the converse sense and, given a family of square-integrable functions, analyze its oscillations at some precise scale that we fix, \eg $h=
2^{-n}$ for $n\in \N$. As illustrated in the next section, this strategy can be used to prove that solutions to dispersive
evolution equations such as wave-type equations or the Schr\"odinger
equation are observable. 

\subsection{High-frequency analysis and  control theory}

On a compact Riemannian manifold $(M, g)$ without boundary,
consider the following  free wave  equation, here of
Klein-Gordon type, 
\begin{align}
  \label{eq: free wave equation}
  \d_t^2 u - \Delta_M u + u =0, 
  \ \ (u,  \d_t u)_{|t=0} = (u_0, u_1). 
\end{align}
It is well-posed for
$(u_0,u_1) \in H^1(M)\times  L^2(M)$. Given an open subset $\omega$
of $M$ and $T>0$, one says that the wave equation is observable from
$\omega$ in time $T>0$
if there exists $C>0$ such that 
\begin{align}
  \label{eq: obs wave equation}
  \mathcal E (u) \leq C \int_0^T \Norm{1_\omega \d_t u }{L^2(M)}^2 dt, 
\end{align}
for any solution $u$ to \eqref{eq: free wave equation}, where
$\mathcal E (u)$ denotes the energy of the solution 
\begin{align*}
  \mathcal E (u)  = \frac12 \big( \Norm{u_0}{H^1(M)}^2 + \Norm{u_1}{L^2(M)}^2 \big).
\end{align*}
With a duality argument \cite{JLL},  an observability inequality as in
\eqref{eq: obs wave equation} is equivalent to the exact
controllability of the wave equation from $\omega$ in time $T$, that
is, for any initial and final states, $(y_{0}, y_{1})$ and $(y^T_0,
y^T_1)$ both in  $H^1(M)\times L^2(M)$, 
the ability to find $f \in L^2((0,T)\times M)$ such that the solution $y$
to 
\begin{align*}
  \d_t^2 y - \Delta_M y +y =1_\omega f, 
  \ \ (y, \d_t y)_{|t=0} = (y_0, y_1), 
\end{align*}
satifisfies $(y, \d_t y)_{|t=T} = (y^T_0,y^T_1)$.

As shown in \cite{Lebeau:92,BDLR:23a}, for the proof of \eqref{eq: obs wave
  equation} it suffices to consider sequences of waves $(u_k)_{k\in\N}$ with
localized time-frequency $\tau \sim h^{-1} \sim 2^n$, $n\in\N$, built by means of the
eigenfunctions $\varphi_k$ defined in~\eqref{eq:eigen}, with
$\sqrt{E_k} \sim h^{-1}$. Although not intrinsic to the considered
question, Section~\ref{sec:app_control} discusses how  a semi-classical point of view can be chosen, offering a powerful analysis toolbox.


\section{Correspondence principle}\label{sec:cor_pri}

The phase space of quantum mechanics is the set~$\R^{2d}$ of positions and momenta: 
\[
z=(x,\xi)\in\R^{2d}.
\] 
The Fourier transform $f\mapsto \widehat f$ is given by 
\[
\widehat f(\xi)= \int_{\R^d} f(x) {\rm e}^{-ix\cdot \xi} dx,\;\;\xi\in\R^d,
\]
and $f \mapsto \mathcal F (f) = (2\pi)^{-\frac d2} \widehat f$ is a
unitary transformation of $L^2(\R^d)$.  In a semi-classical context,
one rescales the Fourier transform by considering the $h$-Fourier
transform $f\mapsto \mathcal F_h f$
\[
\mathcal F_h f(\xi)= (2\pi h)^{-\frac d2}\widehat f\left(\frac \xi h\right),\;\;\xi\in\R^d.
\]

Note that by the Plancherel theorem, the average momentum introduced in~\eqref{def:av_imp} reads 
\[
\langle \xi_j\rangle_{\psi}=\int_{\R^d}\xi_j|\mathcal F_h \psi (\xi)|^2 d\xi.
\]

The phase space  $\R^d\times \R^d$ is endowed with the symplectic form
$\omega=d\xi\wedge dx$ defined by
\begin{equation}\label{def:symplectic_form}
\omega(z,z')= Jz \cdot z',\;\; J=\begin{pmatrix} 0 & {\rm Id}_d \\ -{\rm Id}_d & 0 \end{pmatrix},\;\;z,z'\in\R^{2}.
\end{equation}
Geometrically, it is natural to view the phase space as the cotangent bundle $T^*\R^d$, with $\xi\in T^*_x\R^d$, the cotangent variable (see Section~\ref{sec:app_geo}).

\subsection{Semi-classical wave packets}

Semi-classical wave packets are wave functions associated with a classical state $z= (q,p)\in\R^{2d}$. 
One defines Gaussian wave packets as
\[
g^ h_z(x) = (\pi h)^{-d/4} \exp(-\tfrac{1}{2 h}|x-q|^2 + \tfrac{i}{h}p\cdot(x-q)),
\]
for $x\in\R^d$. It
is normalized, $\Norm{g^ h_z}{L^2}=1$, and centered in $z$, 
\[
\langle x_j\rangle_{g^ h_z} = q_j
\ \ \text{and}  \ \ 
\langle \xi_j \rangle_{g^ h_z} = p_j,\;\; 1\leq j\leq d.
\]
Moreover, its $h$-Fourier transform has the same structure 
\[
  \mathcal F_h\left({\rm e} ^{i\frac {p\cdot q}{2h}} g^h_z\right)
  = {\rm e} ^{i\frac {p\cdot(- q) }{2h}} g^h_{Jz}, 
  \quad z=(q,p)\in\R^{2d}.
\]
\smallskip 

The Gaussian wave packets were introduced in part because they have the unique property
among $L^2$-functions of saturating the uncertainty principle
\[
d_{g^h_z} x_j = d_{g^h_z} \xi_j = \sqrt {\frac h 2}, \;\; 1\leq j\leq d.
\]

Besides, any wave function can be written as a superposition of  Gaussian wave packets according to the  Bargmann formula: for all $f\in L^2(\R^d)$
\begin{equation}\label{eq:bargmann}
f =(2\pi h)^{-\frac d2} \int_{\R^{2d} } \mathcal B_h[f](z) g^h_z dz,
\end{equation}
where the Bargmann transform~\cite{corobook} is the isometry from $L^2(\R^d)$ into   $L^2(\R^{2d})$ defined by
\begin{equation*}
\mathcal B_h[f](z)=(2\pi h)^{-\frac d2} \inp{f}{g^h_z}_{L^2},\;\;z\in\R^{2d}.
\end{equation*}

\subsection{Semi-classical pseudodifferential operators and related notions}

A question that arises from quantum mechanics is the quantization problem, or how to associate an operator to an energy, also called Hamiltonian. It  gives a mathematical setting to  explore the correspondence between classical and quantum mechanics. 

\subsubsection{Quantization of observables }
Let $a(x,\xi)$ be a semi-classical observable in the Schwartz space
$\S(\R^{2d})$. The {\it semi-classical pseudodifferential operator} (\hpdo), of symbol $a$ is the operator $\Op(a)$ defined on functions $f\in \S(\R^d)$ by
\begin{align*}
&\Op(a) f(x)\\
\nonumber
&=(2\pi h)^{-d} \int_{\R^{2d}} a\big(\tfrac12(x+y),\xi \big) {\rm e}^{{i\over h} \xi \cdot (x-y)} f(y) dy\,d\xi.
\end{align*}
This form is called the {\it Weyl-quantization} of the symbol~$a$
\cite{Hor,DimassiSjostrand,MartinezBook,Zwobook}.

The operator $\Op(a)$ maps $\S(\R^d)$ into itself and, by duality, $\S'(\R^d)$
into itself. Its kernel $k_h$ can be expressed in terms of the inverse Fourier transform of~$a$ in the variable~$\xi$
\begin{equation}\label{def:kappa}
\kappa(x,v)= (2\pi)^{-d} \int_{\R^d} a(x,\xi){\rm e}^{i\xi\cdot v}d\xi ,\;\; (x,v)\in\R^{2d}.
\end{equation}
Indeed,  one has 
\[
k_h(x,y) = \frac1{h^{d}}\, \kappa \left( \frac{x+y}2, \frac {x-y} h\right),\;\; (x,y)\in\R^{2d}.
\]
As a consequence of the Schur Lemma, the operator $\Op(a)$ maps $L^2(\R^d)$ into itself and 
\begin{align*}
&\| \Op(a)\|_{\mathcal L(L^2(\R^d))}  \leq  \int_{\R^d}\sup_{x\in\R^d} |\kappa(x,v)| dv\\
&\;\;\leq C \sup_{\substack{\beta\in \N^d\\|\beta|\leq d+1}} \sup_{x\in\R^d} \|\partial_\xi^\beta a(x,\cdot) \|_{L^1(\R^d)},
\end{align*}
for $C>0$ independent of $a$ and $h$.
The Calder\'on-Vaillancourt theorem~\cite{CV,hwang,corobook} also gives the existence of $C>0$ such that for all $a$ and $h$,
\begin{align*}
\nonumber
&
\left\| \Op(a)\right\|_{{\mathcal L}(L^2(\R^d))}\\
&\;\;\leq C
\sum_{\alpha\in\N^{2d},|\alpha| \leq 2d+1} h^{\frac{|\alpha|}2}\sup_{\R^d\times\R^d}|\partial_{x,\xi}^\alpha a|.
\end{align*}
This estimate can be derived from the case $h=1$ by conjugating
$\Op(a)$ by the scaling unitary operator $T_h: f\mapsto h^{\frac d4}
f(\sqrt h\cdot)$. Indeed, one has $T_h \Op(a) T_h^*= {\rm
  Op}_1(a(\sqrt h\cdot,\sqrt h\cdot))$.

\smallskip 
The present definition of {\hpdo}s can be set within the general
H\"ormander formalism with the phase space metric $|dx|^2 + h^2
|d\xi|^2$; see \cite[Sections 18.4-5]{hormander}, \cite[Section 2]{L} and \cite[Sections 2.2--2.3]{MartinezBook}.

\subsubsection{Symbolic calculus}    

The set of {\hpdo}s is an algebra that enjoys symbolic calculus. If $a,b\in\S(\R^{2d})$, then in ${\mathcal L}(L^2(\R^d))$,
\begin{align}
\nonumber 
&\Op(a)\Op(b) =  \Op(ab)\\
\label{eq:produit}
&\;\;+{ h\over 2i}\,\Op \left(\{a,b\}\right)+O\left( h^2 \right),
\end{align}
where $\{a,b\}$ denotes the {\it Poisson bracket}
\begin{align*}
\{ a,b\}=\nabla_\xi a\cdot \nabla_x b- \nabla_x a\cdot \nabla_\xi b.
\end{align*}
This implies that the commutator of two {\hpdo}s
 is of lower order, which turns out to read 
\begin{equation}\label{eq:commutator}
\left[\Op(a),\Op(b) \right] =  \frac{ h}{i}\,\Op (\{a,b\}) +
O\left( h^3\right),
\end{equation}
because of the symmetries of the term  $O(h^2)$ in~\eqref{eq:produit}.

The remainder terms  $O(h^2)$, $O(h^3)$ appearing in~\eqref{eq:produit} and~\eqref{eq:commutator}, involve Schwartz semi-norms of the symbols~$a$ and~$b$, such as
\begin{equation*}
N_{k}(a)= \sup_{|\gamma|\leq k}\| \partial^\gamma_z a \|_{L^\infty}.
\end{equation*} 
for $k\in\N$ large enough~\cite{Robert}.

Regarding the adjoint, one simply has 
\begin{align}
  \label{eq: adjoint}
  \Op(a)^* = \Op(\overline a).
\end{align}
In particular, if $a$ is real-valued, then $\Op(a)$ is a symmetric
bounded operator, thus self-adjoint. 
Results of this section can be found in \cite{DimassiSjostrand,Zwobook,AFF}, for example.
\smallskip

Other quantizations  also enjoy a symbolic calculus. Let us cite the
{\em left-quantization} \cite{Fedoryuk_Maslov}, so-called {\em classical quantization}, $a\mapsto a(x,hD)$ defined by
\begin{align*}
&a(x,hD) f(x)\\
\nonumber
&\quad =(2\pi h )^{-d} \int_{\R^{2d}} a\big(x, \xi \big) {\rm e}^{\frac{i}{h} \xi \cdot
  (x-y)} f(y) dy\,d\xi\\
&\quad =(2\pi )^{-d} \int_{\R^{d}} a\big(x, h\xi \big) 
  {\rm e}^{i \xi \cdot x} \widehat{f}(\xi)\,d\xi, 
  \quad f \in \S(\R^d).
\end{align*}
However, the symbol for the adjoint operator is not as simple as in
\eqref{eq: adjoint} and the remainder in the counterpart to
\eqref{eq:commutator} is only $O(h^2)$ in the left-calculus. This is a
reason for the Weyl-quantization to be often preferred. Correspondance
between the two quatizations is expressed by 
\begin{align*}
  a (x, \xi) = e^{\frac{i h}{2}D_x\cdot D_\xi } b (x,\xi),
\end{align*}
if $a(x, hD) = \Op(b)$,~\cite{DimassiSjostrand}.
\smallskip 

The notations $a(x,hD)$ and $\Op(a)$ are extended to smooth functions $(x,\xi)\mapsto a(x,\xi)$ that satisfy symbol estimates of the form 
\begin{equation*}
\forall \alpha,\beta\in\N^d,\;\exists C_{\alpha,\beta}>0,\;\left\| \langle \xi\rangle^{-m+|\beta|} \partial_x ^\alpha \partial_\xi^\beta a\right\|_{L^\infty} \leq C_{\alpha,\beta}
\end{equation*}
for some $m\in\N$ (here $\langle \xi\rangle= \sqrt{1+| \xi|^2}$). One then says that $a\in S^m$~\cite{Zwobook}. 

In particular, this class contains  the
functions~$p$ that are polynomial functions of degree~$m$ in the
variable~$\xi$ with coefficients that are smooth bounded functions
of~$x$ with bounded derivatives. In this case, the operators $p(x,hD)$ and
$\Op(p)$  are differential operators. 

For such symbol classes, symbolic calculus results above also hold. In
particular, 
if $a \in S^m$ and $b \in S^{m'}$,
then $\Op(a) \Op(b) = \Op(c)$ with $c \in S^{m+m'}$ given by $c = ab +
h \{ a,b\} /(2i) \mod  h^2 S^{m+m'-2}$. 

Introducing the  semi-classical Sobolev norms
\[
\| f\|_{s}= \sup_{0\leq \ell \leq s} \| \langle hD_x\rangle ^\ell f\|_{L^2},\;\; s\in \R,
\]
if $a \in S^m$ and $s\in \R$, 
there exists a constant $C>0$ such that 
\[ 
\| \Op(a)f \|_{s} \leq C\| f\|_{s+m},\quad h\in(0,1], \  f\in\S(\R^d).
\]

\subsubsection{Bargmann transform and {\hpdo}s}

The relations of~{\hpdo}s   with the Bargmann transform enlighten the
role of the {\hpdo}s in terms of microlocalization. For $a \in
\S(\R^{2d})$,  there exists a constant $C>0$ such that for $h\in(0,1]$, 
\begin{align}\label{link_bargmann}
\left\| \Op(a)-\mathcal B_h^* a\mathcal B_h\right\|_{\mathcal L(L^2(\R^d))}\leq Ch.
\end{align}
Indeed, the kernel of the operator $\mathcal B_h^* a\mathcal B_h$ is the function 
\[
k_h^{\mathcal B}(x,y)=  \frac1{h^{d}}\, \kappa_h^{\mathcal B} \left( \frac{x+y}2, \frac {x-y} h\right),\;\; (x,y)\in\R^{2d},
\]
 related with the function $\kappa$ of~\eqref{def:kappa}  according to 
\[
\kappa_h^{\mathcal B}(x,v)=\pi^{-\frac d2}{\rm e}^{-\frac h4|v|^2} \int_{\R^d} \kappa(x- \sqrt h q,v){\rm e}^{-|q|^2} dq,
\]
for  $(x,v)\in\R^{2d}$.
Therefore, using Taylor expansions, the fact that $\int q e^{-|q|^2}
=0$, and the rapid decay of $\kappa(x,v)$ in $v$ one obtains 
\[
\kappa_h^{\mathcal B}(x,v) -\kappa(x,v)=h\int_{\R^d} A_h(x,q,v){\rm e}^{-|q|^2} dq,
\]
where for all $N\in\N$, the function 
\[
(x,q,v)\mapsto |v|^N A_h(x,q,v)
\]
is uniformly bounded in $h\in(0,1]$. Estimate~\eqref{link_bargmann} then comes from the Schur Lemma.

\subsubsection{Ellipticity, parametrix, and sharp G{\aa}rding inequality}

Symbolic calculus allows one to transfer
 properties of the symbol $a$ to the  $h$-$\psi$do  $\Op(a)$.

Let $P^h=p(x,hD)$ be a  differential operator with a symbol $p(x,\xi)$
that is a smooth polynomial function of degree~$m$ in $\xi$
\[
p(x,\xi)=\sum_{|\alpha|\leq m} p_\alpha(x) \xi^\alpha.
\]
One has $p \in S^m$. 
The symbol $p$ is said to be {\it elliptic} if there exists $C>0$ and
$R>0$ such that 
\begin{align*}
  |p(x,\xi)| \geq C |\xi|^m, \ (x, \xi) \in \R^{2d}, \ |\xi| \geq R. 
\end{align*}
\smallskip 
In such a case, the \hpdo~$P^h$ is one to one  from $H^{s+m}_h(\R^d)$ onto $H^s_h(\R^d)$ for all $s\in \R$ and 
\[
(P^h)^{-1}= \Op(p^{-1})+O(h), 
\]
by symbolic calculus. One has $p^{-1}\in S^{-m}$ 
and $\Op(p^{-1})$ is called a {\it parametrix} of $P^h$. 

The question of positivity is addressed by the sharp {\it G\aa rding inequality}, which is a direct consequence of estimate~\eqref{link_bargmann}.
There exist $C,N>0$ such that  for all $a$ in $\S(\R^d)$ satisfying $a\geq 0$, we have for all $f$ in $\S(\R^d)$ and $h$ in $(0,1]$.
\begin{align}
\label{eq: gaarding}
\left( f,\Op(a)f\right)\geq -Ch \|f\|_{L^2} \sup_{|\alpha|\leq N} \|\partial^\alpha_z a\|_{L^\infty}.
\end{align}

\subsubsection{Functional calculus and trace formula}

Since $\Op(a)$ is a bounded self-adjoint operator for real-valued $a$ in $\S(\R^d)$,
 functional calculus can be used and the operator $F(\Op(a))$ is
well defined for $F$ continuous on $\R$.
\smallskip 

 Suppose $F\in\Cinfc(\R)$. Then,
$F(\Op(a))$ coincides asymptotically with a pseudodifferential
operator of symbol $F(a)$, that is, 
\begin{equation}\label{Funct_calc}
F(\Op(a))=\Op(F(a)) +O(h) \ \ 
\text{in} \ \ {\mathcal L}(L^2(\R^d)).
\end{equation}

This relies on the {\it Helffer-Sj\"ostrand formula}~\cite{Zwobook,DimassiSjostrand} that plays an
important role in semi-classical analysis and is of interest in
itself, in particular because of the alternative  construction of the
functional calculus it provides for a (possibly unbounded)
self-adjoint operator~\cite{Davies}.
\smallskip 

In fact, for all  $n\! \in\! \N$, $F$ has an almost analytic
continuation, that is,  a function $\tF_{n}\in \Cinfc(\C)$
that coincides with $F$ on $\R$ and such that 
\begin{equation}
  \label{eq: almost-analytic}
  \left|  {\bar{\d}}
    \tF_{n}(z)\right|
  \leq C\, |\Im (z)|^{n}, \quad z\in \C.
\end{equation}
The Helffer-Sj\"ostrand formula reads
\begin{equation*}
  F\big(\Op(a)\big)
  ={1\over\pi} \int_{\C}\bar{\d}
  \tF_{n}(z)\big(\Op(a)-z\big)^{-1}
L(dz),
\end{equation*}
where  $L(dz)$ is the Lebesgue measure on $\mathbb C$.  The operator $\big(\Op(a)-z\big)^{-1}$ is bounded, with norm $|\Im(z)|^{-1}$ for almost all $z\in\C$ and, 
thanks to \eqref{eq: almost-analytic}, using a parametrix of $\Op(a)-z$ to replace the  resolvent
$(\Op(a)-z)^{-1}$ one obtains~\eqref{Funct_calc}.
\smallskip

Noticing that  for all fixed $h>0$, $\Op(a)$  is a compact operator with Hilbert-Schmidt norm
\[
\| \Op(a)\|_{{\rm HS}(L^2(\R^d))} = (2\pi h)^{-d/2} \|a\|_{L^2(\R^{2d})},
\]
one deduces a trace formula: for $F\!\in\! \Cinfc(\R)$ nonnegative one has
\begin{align}\label{trace}
& \mathop{\rm Tr} \big(F(\Op(a)\big) \\
\nonumber
&\qquad \mathop{\sim}_{h\rightarrow 0} \, (2\pi h)^{-d} \int_{\R^{2d}} F(a(x,\xi)) dx d\xi.
\end{align}

This approach is used in the spectral analysis of Schrödinger
operators such as $-h^2 \Delta +V(x)$ for confining potential, or magnetic Schrödinger operators  $-|hD_x-A(x)|^2$ on bounded domains (see the historical series of papers by B.~Helffer and J.~Sj\"ostrand~\cite{HS1,HS2,HS3,HS4} and the books~\cite{Fournais_Helffer,Raymond_book,SanVuNgoc_book}).

\subsection{Wigner transform and semi-classical measures }
\subsubsection{Main definitions and example}
Following E.~Wigner~\cite{Wigner}, once given a bounded family $(\psi^h)_{h>0}$ in $L^2(\R^d)$, one can
consider the distribution 
\[
 W[\psi^h]: a\mapsto \dup{W[\psi^h]}{a}= (\Op(a)\psi^h,\psi^h)
\]
called the {\it Wigner transform} of  $(\psi^h)_{h>0}$. One finds it
is defined for $(x,\xi)\in\R^{2d}$ by 
\begin{align}\label{def:wigner_transform}
&W [\psi^h](x,\xi) =(2\pi)^{-d} \int_{\R^d}{\rm e}^{iv\cdot \xi}  \\
\nonumber
&\qquad \times\psi^h\left(x-\frac h 2 v\right) \overline \psi^h\left( x+\frac h 2 v\right) dv.
\end{align}
This notion has been revisited in the 1990's, see~\cite{HMR87,LionsPaul} and the works of P. Gérard and his coauthors~\cite{Ge91,GerLeich93,GMMP}.
\smallskip

In view of~\eqref{eq: gaarding}, for any family, weak limits point in
the sense of distributions of the Wigner transform of
$(\psi^h)_{h>0}$ are  finite nonnegative
measures. They
are called {\it semi-classical measures} of the family~$(\psi^h)_{h>0}$ (see~\cite{HMR87,Ge91,GerLeich93}). One also uses the term {\it Wigner measures} (see~\cite{GMMP}).  For such a measure $\mu$,  there exists a subsequence $h_k\Tend{k}{+\infty}0$ such that 
\begin{equation}\label{eq:scm}
\dup{W[\psi^{h_k}]}{a}
\Tend{k}{+\infty} \dup{\mu}{a},\;\;\forall a\in \Cinfc(\R^{2d}). 
\end{equation}

For example, the Wigner transform of the Gaussian wave packet $g^h_z$ is given  for $z,\zeta\in\R^{2d}$ by
\[
W[g^{h}_z](\zeta) = (\pi h)^{-d} \exp(-\tfrac{1}{h}|\zeta-z|^2).
\]
Thus, the family $(g^h_z)_{h>0}$ has only one semi-classical measure,
namely, 
\[
\mu(x,\xi)= \delta(x-q)\otimes \delta(\xi-p).
\]
In the limit $h\rightarrow 0$, the wave function $g^h_z$ converges to the classical state $z=(q,p)$, which gives a first illustration of the correspondence principle.

\subsubsection{$h$-oscillation}

There is a connexion between  the weak limits of $|\psi^h(x)|^2 dx$ and the semi-classical measures of $(\psi^h)_{h>0}$.
Indeed, if the sequence $(h_k)_{k\in\N}$ and the measure $\mu$ fulfills property~\eqref{eq:scm} and if $\nu$ is a weak limit of the measure $|\psi^{h_k}(x)|^2 dx$, then 
\[
\nu\big (\{ x\}\big) \geq 
\mu\big(\{x\}\times \R^d\big)
\]
as measures on $\R^d_x$. 
Besides, equality holds if $(\psi^h)_{h>0}$ is 
$h$-oscillating, namely satisfies the property 
\[
\limsup_{h\rightarrow 0} \int_{h|\xi|\geq R} |\widehat \psi^h(\xi)|^2 d\xi \Tend{R}{+\infty} 0.
\]
In other words, no mass escapes to infinity in frequency.
Such a property is satisfied for examples if  $(\langle hD_x\rangle  ^s \psi^h)_{h>0}$ is uniformly bounded in $L^2(\R^d)$ for some $s>0$. 
In fact, once given a bounded family in
$L^2(\R^d)$, an appropriate semi-classical scale (if any) can be
sought by analyzing
the size of one of its Sobolev norms, motivating a  semi-classical
analysis at that precise scale. Such strategies will be implemented in Sections~\ref{sec:applications} for 
the analysis of the examples  
presented in   Section~\ref{sec: Some semi-classical problems}.

\subsubsection{Wave front set}

The support of the semi-classical measure of a bounded family $(\psi^h)_{h>0}$  in $L^2(\R^d)$ is included in the {\it semi-classical wave front set} denoted ${\rm WF}_h(\psi^h)$. The latter is characterized by the following property: $(x,\xi)\notin {\rm WF}_h(\psi^h)$ if and only if there exists an open  neighborhood $U$ of the point $(x,\xi)$ and a function $a\in\Cinfc(U)$ such that 
\begin{align*}
&a(x,\xi)\not=0\\
&\mbox{and}\;\;
\forall n\in\N , \;\; \| \Op(a) \psi^h\|_{L^2}=O(h^n).
\end{align*}
If $\mu$ is a semi-classical measure of $(\psi^h)_{h>0}$ for the scale~$h_k$, 
\[
{\rm Supp}\, \mu \subset  {\rm WF}_{h_k}(\psi^{h_k}).
\]
Historically, the semi-classical wave front set was introduced earlier than semi-classical measures. It is closely related to microlocal versions of wave front set where no  scale is emphasized (see~\cite[Vol. 1, Ch. 8]{Hoermander:V1}).

\subsubsection{Semi-classical measures and PDEs}
\label{sec: semi-classical measures and PDEs}
Consider $P^h = p(x,h D)$ a differential operator.  Suppose $(\psi^h)_{h>0}$ is a sequence of
bounded $L^2$-functions associated with a semi-classical measure
$\mu$ such that 
\[
P ^h\psi^h = o(1)
\]
 in $L^2(\R^d)$ as $ h \to 0$. Then, for $a\in\S(\R^{2d})$,
\[
\inp{\Op(a) P^h \psi^h}{\psi^h}_{L^2} = o(1),
\]
 implying  
$\dup{\mu}{ap}=0$ and 
\begin{align}
  \label{eq: localization supp mu}
  \supp(\mu) \subset \Char(P^h), 
\end{align}
where $\Char(P^h) = \{ p(x,\xi) =0\}$ is the {\it characteristic set} of $p$. 
\smallskip 

Assume moreover that $P^h$ is symmetric and 
\[
P^h \psi^h = o(h)
\]
 in
$L^2(\R^d)$ as $ h \to 0$. Then, for $a\in\S(\R^{2d})$,
\[
\inp{[\Op(a), P^h] \psi^h}{\psi^h}_{L^2} = o(h),
\] implying
$\dup{\mu}{ \{p,a\} }=0$. 
One has $\{p,a\} = \Hp a$ with $\Hp = J\, 
\nabla_{\! \! x,\xi} \, p$, the
Hamiltonian vector field associated with $p$ (recall that $J$ is given by~\eqref{def:symplectic_form}).  Since $\transp{\Hp} = -\Hp$
one finds
\begin{align}
  \label{eq: transport mu}
  \Hp \mu = 0,
\end{align}
in the sense of distributions, 
meaning with \eqref{eq: localization supp mu} that $\mu$ is
invariant along the Hamiltonian curves $(\Phi^t(z))_{t\in\R}$ for $z\in \R^{2d}$ where the map
$\Phi^t:\R^{2d}\to\R^{2d}$, 
 is determined by
\begin{equation}\label{def:Hp}
\dot\Phi^t = \Hp (\Phi^t),\quad \Phi^0 =\mathrm{Id}_{\R^{2d}}.
\end{equation}
For all $t\in\R$, $z\mapsto \Phi^t(z)$ is a symplectomorphism (it
preserves the symplectic form $\omega$ given
in~\eqref{def:symplectic_form}).  Condition~\eqref{eq: transport
  mu} relates classical phase-space trajectories and solutions
concentrations. 

Propagation of semi-classical measures is more difficult to prove if
coefficients are singular. We refer for instance to \cite{Gannot-Wunsch23,Galkowski-Wunsch23,BDLR:23b}.

\subsection{Semi-classical evolution}

Consider an evolution equation involving a smooth time-dependent Hamiltonian function $p: \R\times \R^{2d}\to\R$,  with sub-quadratic growth
\begin{equation*}
\forall N\geq 2,\;\exists C_N>0,\; \sup_{|\beta|=N} \, \sup_{(t,z)\in\R\times \R^{2d}}\left|  \partial_{z}^\beta p(t,z)\right| \leq C_N.
\end{equation*}
Then,  the operator 
$P^h(t)=\Op(p(t))$ 
 is   self-adjoint   and  there exists a  strongly continuous two-parameters family of unitary operators  
 $U^h (t,s)$ such that
 \begin{equation*}
   i h {d\over dt}U^ h (t,s) =P^h(t) \, U^ h  (t,s),
   \quad U^ h  (s,s)={\rm Id}_{L^2}
\end{equation*}
on the domain of the operator $P^h(t)$ (see \cite{ReedSimon}).  If
$P^h$ is independent of time, then $U^h(t,s) =U^h(t-s,0)$, and
$U^h(t,0)$ is the semigroup generated by $P^h$.

\subsubsection{The Egorov Theorem}

At the classical level, 
one associates with $p(t)$ the ordinary differential system $\partial_t z =\Hamiltonian_{p(t)}(t,z)$ and the flow map $
\Phi^{t,s}:\R^{2d}\to\R^{2d}$, 
that is determined by
\begin{equation}\label{def:Phits}
\partial_t\Phi^{t,s} = \Hamiltonian_{p(t)} (t,\Phi^{t,s}),\quad \Phi^{s,s} =\mathrm{Id}_{\R^{2d}}.
\end{equation}
Note that if $p=p(z)$ does not depend on the time, $\Phi^{t,s}=\Phi^{t-s}$ defined in~\eqref{def:Hp}.

For the evolution of an observable $a\in\S(\R^{2d})$ one uses the Liouvillian $\mathcal L_{t,s}a = a\circ\Phi^{t,s}$ that satisfies the transport equation
\[
\partial_t(\mathcal L_{t,s} a) = \{p(t),\mathcal L_{t,s} a\},\quad \mathcal L_{s,s} a = a.
\]

  At the quantum level, one works with the quantization of $p(t)$, the
  operator $P^h(t)$ and, given an observable~$a$, one considers the
  conjugation of the operator $\Op(a)$ by the propagators $U^h(t,s)$:
\[
U^ h (s,t) \circ \Op(a) \circ U^ h  (t,s).
\]
The Egorov Theorem connects 
the classical picture and the quantum one in the limit $h\rightarrow 0$ (see for instance~\cite{Rb})
 \begin{thm}\label{th:egorov}
There exists a constant $C>0$ such that for all $a\in \S(\R^{2d})$ and $t,s\in \R$
\begin{align}
\nonumber
&\left\| U^ h (s,t) \circ \Op (a) \circ U^ h  (t,s) - \Op ({\mathcal L}_{t,s}a )\right\|_{\mathcal L(L^2(\R^d))} \\
\label{eq:egorov}
&\;\leq C \,  h^2 \, |t-s| \,{\rm e}^{C|t-s|}   N(a),
\end{align}
where $N(a)$ denotes a fixed  semi-norm of $a$. 
\end{thm}

For some $\delta\in (0,1)$, on a large time interval of size $|t-s|\sim \frac  \delta C  \ln
\left(\frac 1 h\right)$, the error estimate in \eqref{eq:egorov} is $\delta h^{2-\delta}
\ln \left(\frac 1 h\right) N(a) \ll 1$. This large time for
which $\Op ({\mathcal L}_{t,s}a )$ provides a good approximation is
called the {\it Ehrenfest time} and characterizes the range of validity of the semi-classical approximation~\cite{Bouzouina_Robert}. 
\smallskip

\subsubsection{Semi-classical measures and propagators}

Assume in this section that the Hamiltonian $p$ does not depend on the time, $p=p(x,\xi)$. To study
\[
\psi^h(t):= U^{h}\left(t,0\right)\psi^h_0,\;\; h>0,
\]
on (possibly large) time scales $t\sim 1/h^{\alpha}$, $\alpha\geq 0$, one considers the limit as $h$ goes to $0$ of the quantities 
\[
\int_{\R} \theta(t) \langle W[\psi^h(t/h^\alpha)] , a\rangle\, dt 
\]
for $\theta\in L^1(\R)$ and $a\in\S(\R^{2d})$. 
Up to the extraction of a subsequence, this limit is described by a family of
measures $d\mu^t_\alpha (x,\xi) \otimes dt$ that is also called a semi-classical measure of the family $(\psi^h(t/h^\alpha))_{h>0}$. 
\smallskip

The Egorov Theorem implies the following:
\begin{description}
\item[Case $\bld{\alpha=0}$.] Any  semi-classical measure $d\mu^t_0 (x,\xi) \otimes dt$ of $\left(\psi^h(t)\right)_{h>0}$ satisfies $\mu^t=\Phi^{t,0}_*\mu$ for $\mu$ a semi-classical measure of $(\psi^h_0)$.
\item[Case $\bld{\alpha>0}$.] Any  semi-classical measure $d\mu^t_\alpha (x,\xi) \otimes dt$ of  $\left(\psi^h\left(\frac t{h^\alpha}\right)\right)_{h>0}$ satisfies the invariance property:
$\mu^t_\alpha=\Phi^{s,0}_*\mu^t_\alpha$ for all $s\in\R$. In other words, the measure 
$\mu^t_\alpha$ is invariant by the flow  $s \mapsto \Phi^{s,0}$.
\end{description} 

When $\alpha=0$, the description of measure given above in this case opens algorithmic strategies for a numerical computation of the Wigner transform of $\left(\psi^h(t)\right)_{h>0}$. At leading order, this Wigner transform is approximated by the Wigner measure, and thus  by the pull-back by the flow $\Phi^{t,0}$ of the Wigner transform of $(\psi^h_0)_{h>0}$, that can be computed numerically via a quadrature procedure for the integral~\eqref{def:wigner_transform}.
 The correspondence principle allows to trade the resolution of a $h$-dependent PDE by solving $h$-independent ODEs \cite{LasserKeller}. 
 \smallskip

In the case  $\alpha>0$, the invariance property 
of $\mu^t_\alpha$ 
   implies that $\supp(\mu^t_\alpha)$ is a union
of periodic orbits of the flow. For example, if   $p = |\xi|^2/2$, the flow $\Phi^{s,0}$ is
 given by $(x,\xi)\mapsto (x+s\xi,\xi)$; the fact that the measure $\mu^t_\alpha$ is of finite mass and invariant by $\Phi^{s,0}$ implies
 $\supp(\mu^t_\alpha) \subset \{\xi=0\}$; this  illustrates the dispersion
 effects in the Schr\"odinger equation. 
 Such an analysis is at the roots of the results of~\cite{AM:14} on the torus, for example.

\subsubsection{Propagation of coherent states}

The propagation of coherent states also illustrates the correspondence principle.
For
$z=(q,p)\in \R^{2d}$, the function $U^h(t,s) g^h_z$ can be described
at leading order via classical quantities. We need to introduce additional notations
\smallskip 

One enlarges the set of profiles  and considers complex-valued Gaussian profiles~$g^\Gamma$, whose  
covariance matrix~$\Gamma$ is taken in the Siegel half-space ${\mathfrak S}^ +(d)$ of  $d\times d$ complex-valued symmetric matrices with positive imaginary part,
\[
{\mathfrak S}^+(d) = \left\{\Gamma\in\C^{d\times d},\
  \Gamma=\transp\, \Gamma,\ \Im\Gamma >0\right\}.
\]
More precisely, $g^\Gamma$ is given by
\begin{equation*}
g^\Gamma(x)
:= c_\Gamma\, {\rm e}^{\frac{i}{2}\Gamma x\cdot x},\quad x\in\R^d,
\ \Gamma\in{\mathfrak S}^+(d),
\end{equation*}
where 
$c_\Gamma=\pi^{-d/4} {\rm det}^{1/4}(\Im\Gamma)$  
is a $L^2$-normalization constant. 

For $z=(q,p)\in\R^{2d}$, set
 \[
 g^{\Gamma,h}_z(x)=h^{-\frac d4} {\rm e}^{\frac ih p\cdot (x-q)}
 g^\Gamma\left(\frac{x-q}{\sqrt h}\right),\;\; x\in\R^d.
 \]
 Note that $g^{i{\rm Id}_d,h}_z=g^h_z$.
 \smallskip 
 
We also  introduce classical quantities associated with
the  flow map $\Phi^{t,s}$ introduced in~\eqref{def:Phits}.
 Firstly, consider the $d \times d$ blocks of  {\it the Jacobean matrix}
$F(t,s,z) = \partial_{z}\Phi^{t,s}(z)$ 
\begin{align*}
F(t,s,z) = \begin{pmatrix}  A(t,s,z) &B(t,s,z) \\ C(t,s,z) &D(t,s,z)\end{pmatrix},
\end{align*}
which satisfies the linearized flow equation
\begin{align*}\label{eq:F}
\partial_t F(t,s,z) &= J {\rm Hess}_z p (t,\Phi^{t,s}(z)) \, F(t,s,z),
\end{align*}
with 
$F(s,s,z) = {\rm Id}_{2d}.$
The matrix-valued function $F$ is smooth in $t,s,z$ with any derivative in $z$ bounded. 

Secondly, we introduce {\it the action integral}
\begin{equation*}
S(t,s,z) = \int_{s}^t \left(\xi(t')\cdot \dot x(t')-p(t',z(t'))\right)dt',
\end{equation*} 
where we have set $z(t)=\big(x(t),\xi(t)\big) = \Phi^{t,s}(z)$. 
\smallskip

With this notation,  for $\Gamma\in \mathfrak S^+(d)$, one has in~$L^2(\R^d)$
\begin{equation}\label{wp_app}
U^h(t,s) g^{\Gamma, h}_{z}= {\rm e}^{\frac i h S(t,s,z)} g^{\Gamma(t,s,z), h}_{\Phi^{t,s}(z)} + O(\sqrt  h ),
\end{equation}
 with  
\begin{align*}
\Gamma(t,s,z) 
  &= (C(t,s,z)+ D(t,s,z)\Gamma)\\
    &\quad \quad \times (A(t,s,z) +B(t,s,z)\Gamma)^{-1}.
\end{align*}
Having $\Gamma(t,s,z) \in S^+(d)$ follows from (non elementary)
algebraic relations.  The description can be made more precise
with an asymptotic expansion in powers of
$\sqrt h$~\cite{corobook}.  \smallskip

The propagation of semi-classical wave packets was also investigated
in nonlinear contexts by various authors.  
Wave packets are flexible enough for some nonlinear
superposition results to hold. We refer to the book of
R.~Carles~\cite{C} and the references therein.

\subsubsection{Semi-classical approximation of the propagator}

The description of the propagation of Gaussian states and the formula~\eqref{eq:bargmann} yield approximation formulae for the propagator that can be used for a numerical determination of $U^h(t,0)\psi$, $\psi\in L^2(\R^d)$.

One defines the action of the  {\it thawed Gaussian approximation} on $\psi\in L^2(\R^d)$ by  
\begin{align*}
{\mathcal I}^h_{\rm th}(t)\psi =(2\pi h)^{-d} \int_{\R^{2d}} 
\langle \psi, g^ {h}_z\rangle   {\rm e}^{\frac{i}{ h}S(t,0,z)} g^ {h,\Gamma(t,0,z)}_{\Phi^{t,0}(z)} dz ,
\end{align*}
and the {\it frozen Gaussian approximation} by
\begin{align*}
&{\mathcal I}^h_{\rm fr}(t)\psi\\
&\quad =  (2\pi h)^{-d} \int_{\R^{2d}} 
\langle \psi, g^ h_z\rangle  k(t,0,z) {\rm e}^{\frac{i}{ h}S(t,0,z)} g^ h_{\Phi^{t,0}(z)} dz,
\end{align*}
with 
\begin{multline*}
k(t,0,z) = 2^{-d/2}  {\rm det}^{1/2} \Bigl( A(t,0,z) +D(t,0,z) 
\\+i(C(t,0,z) -B(t,0,z) )     \Bigr),
\end{multline*}
which has the branch of the square root determined by continuity in time. The operator 
${\mathcal I}^h_{\rm fr}(t)$ is often referred to as the {\it
  Herman-Kluk propagator}, see~\cite{RS,R,Kluk}. 
\smallskip 

The operators ${\mathcal I}^{ h}_{\rm th/fr}(t) $, built on classical quantities, approximate   the unitary propagator $U^ h(t,0)$, 
giving another illustration of the correspondence principle.

\begin{thm}[\cite{RS,R}]\label{thm:scalarhk}
Let $p(t)$ be a smooth sub-quadratic Hamiltonian, then for $h\in(0,1]$,
\[
  U^ h (t,0) = {\mathcal I}^{ h}_{\rm th/fr}(t) +  O( h)
  \quad \text{in} \ \mathcal L\big( L^2(\R^d)\big).
\]
\end{thm} 

This result illustrates one of 
the paradigms of the semi-classical approach, consisting in trading
the resolution of oscillating PDEs for that of
ODEs.

Note that the thawed/frozen Gaussian operators are order $h$
approximation of the propagator while the  wave packet approximation
of~\eqref{wp_app} is of order $\sqrt h$. This comes from the structure
of the remainder term in~\eqref{wp_app} and integration in $z$. A numerical implementation of
this approximation was carried out in~\cite{LS}. 
\smallskip 

The operators ${\mathcal I}^{ h}_{\rm th/fr}(t)$ belong to the class
of  {\it Fourier integral operators} (FIO). Designing operators {that approximate} the dynamics of a semi-classical propagator  goes back to the early days of semi-classical analysis, see J.~Chazarain \cite{Ch}, B.~Helffer and D.~Robert \cite{HeRo} and \cite{Rb}, see also the books~\cite[Chapter 12]{Zwobook} or~\cite{DimassiSjostrand}.

\section{Applications}\label{sec:applications}

\subsection{Semi-classical analysis of molecular dynamics }\label{sec:app_BO}

\subsubsection{Square integrable families valued in Hilbert spaces}

The semi-classical pseudodifferential calculus naturally extends to the space $L^2(\R^d, \mathcal H)$ for  some Hilbert space~$\mathcal H$, such as $\C^N$ or $L^2(\mathbb T^d)$ where $\mathbb T^d$ is the $d$-dimensional torus, for example. One then proceeds as follows:
\begin{enumerate}
\item [(i)] The symbols $a$ are smooth compactly supported functions from $\R^d$ into the set $\mathcal K(\mathcal H)$ of compact operators on $\mathcal H$,
\item[(ii)] The semi-classical measures are characterized by a positive measure~$\mu$ and a measurable family $M$ defined on $\R^{2d}$ and valued in the set of operators on $\mathcal H$ that are $d\mu$-a.e. nonnegative  trace-class operators~\cite{GerardMDM91}. 
\end{enumerate}

Then, if $(\psi^h)_{h>0}$ is uniformly bounded in $L^2(\R^d,\mathcal H)$, the pair $(M,\mu)$ is a semi-classical measure of $(\psi^h)_{h>0}$ if, up to a subsequence,  for all $a\in\Cinfc(\R^d, \mathcal K(\mathcal H))$, 
\[
(\Op(a)\psi^h,\psi^h)\Tend  {h}{0} \,\left\langle {\rm Tr}_{\mathcal L(\mathcal H)} (a M ),\mu\right\rangle.
\]

Taking $\mathcal H=L^2(\mathbb T^d)$ turns out to be pertinent for the study of periodic problems (see~\cite{CFM3}).
Taking
$\mathcal H=\C^N$ leads to   the  framework of  the Schr\"odinger equation~\eqref{eq:schro_BO} in the Born-Oppenheimer approximation. The symbols then are matrix-valued and the semi-classical measures are characterized by Hermitian matrices~\cite{GMMP}.

\subsubsection{Molecular dynamics}

Consider $U^h(t)$ the unitary propagator associated with equation~\eqref{eq:schro_BO}.
Denote by ${\rm sp}V(x)$ the set of the eigenvalues of the self-adjoint matrix $V(x)$ and let
$\lambda(x)$ be an eigenvalue of $V(x)$ such that 
\begin{align}\label{def:adiabatic}
&\exists \delta_0>0,\; \forall x\in\R^d,\\
\nonumber 
&\qquad \qquad \;{\rm dist}\left(\lambda(x), {\rm sp}V(x)\setminus\{\lambda(x)\}\right)>\delta_0.
\end{align}
Denote by $\Pi(x)$ the associated (smooth) eigenprojector:
\[
V(x) \Pi(x)=\Pi(x) V(x)=\lambda(x) \Pi(x),\;\;\forall x\in\R^d.
\]
Denote by $\Phi^t$ the classical flow associated with the scalar Hamiltonian
\[
p(x,\xi)= \frac {|\xi|^2} 2 + \lambda(x),
\]
as in~\eqref{def:Hp},
and denote by $\mathcal L^t$ the associated Liouvillian, $\mathcal L^t : a\mapsto a\circ \Phi^t$.
Matrix-valued aspects are treated by  introducing the  {\it parallel transport} of matrices along the flow. Let 
\[
F(x,\xi):= 
\left[ \xi\cdot\nabla \Pi(x)\, ,\, \Pi(x)\right],
\] 
  and consider the unitary transforms $\mathcal R(t,z)$ defined for
  $t\in\R$, $z\in\R^{2d}$ by
  \[
  \partial_t \mathcal R(t,z) = F \left( \Phi^{t}(z) \right) \mathcal R(t,z),\;\; \mathcal R(0,z)={\rm Id}.
  \]
  The map $\mathcal R (t,z)$ preserves the eigenspaces along the flow and maps a vector $\vec V_0$  which is in the range of $\Pi_j(x_0)$ to a 
vector ${\mathcal R}(t,z_0) \vec V_0$  in the range of $\Pi_j(\Phi^{t}_j(z_0))$, $z_0=(x_0,\xi_0)$.
 \smallskip 
 
With these notations in hand, the Egorov Theorem admits the following extension~\cite{corobook} to adiabatic situations.
\begin{thm}\label{thm:adiabatic}
Assume~\eqref{def:adiabatic}, 
 then there exists a constant $C>0$ such that for all $a\in \S (\R^{2d},\C^{N,N})$ and $\theta\in L^1( \R)$
\begin{multline}\label{adiabatic}
\Bigl\| \int_{\R} \theta(t) \Bigr(U^ h (-t) \circ \Op (\Pi a\Pi) \circ U^ h  (t) \\
- \Op (\Pi {\mathcal L}^{t} (\mathcal R(-t) \,a \,\mathcal R (-t)^*)\Pi )\Bigr) dt \Bigr\|_{\mathcal L(L^2(\R^d))} 
\\
\leq C \,  h \,  N( a),
\end{multline}
where $N(a)$ denotes a fixed  semi-norm of $a$.

\end{thm}
Property~\eqref{def:adiabatic} is called {\it adiabaticity}, from the
greek a-diabatos = impassable, because, at leading order, the
propagation holds inside the eigenmode defining the Hamiltonian
$p(x,\xi)$. Note also that the generalization of the Egorov theorem requires averaging in time. 
 \smallskip 
 
This result
extends to general time dependent subquadratic Hamiltonians
$p=p(t,x,\xi)$ with eigenvalues  and eigenprojectors that depend
simultaneously on the position and momentum variables
up to  the
introduction of classical quantities associated to each
eigenvalues
 \cite{corobook}.
 \smallskip 

The proof of~\eqref{adiabatic} relies on a diagonalization process using what is called {\it super-adiabatic projectors}, as carried out by A.~Martinez and V.~Sordoni~\cite{MS} {as well as } H.~Spohn and S.~Teufel~\cite{ST}, see G.~Nenciu's work~\cite{N1} for earlier results. See also~\cite{sjo_adiab}.

As for scalar equations, one can extend the  thawed/frozen Gaussian approximations and construct FIO approximating the propagator $U^h(t)$ associated with $P^h$ by using the eigenprojector and the classical quantities associated with the eigenvalues~\cite{FLR}.  Let us also mention nonlinear results for systems in~\cite{C,H} and for initial data that are semi-classical wave packets (see also references therein).

\subsubsection{Eigenvalue crossings}

If the adiabatic condition~\eqref{def:adiabatic} is not satisfied, or
if the gap between the eigenvalues shrinks as $h\to 0$ (see \cite{HJ2}), transitions
between modes may occur. These non adiabatic effects were observed in
the early 1930s by L.~Landau~\cite{La} and C.~Zener~\cite{Ze} independently.
They were investigated more in details in the 1990's, 
starting with the work of G.~Hagedorn~\cite{Hag94} for the
equation~\eqref{eq:schro_BO} with Gaussian wave
packets for  initial data.  \smallskip

For a general  Hamiltonian $H(t,x,\xi)$, crossings were classified in
the early 2000's by Y.~Colin de Verdi\`ere~\cite{CdV} through a
reduction to normal forms. The analysis of the semi-classical measures
and Wigner transforms is understood in these generic
situations~\cite{FG2,F04,FL1}.  The loss of adiabaticity led to
replace the Liouville operator of Theorem~\ref{thm:adiabatic}, by a
Markov process including branches of classical trajectories and a
branching procedure whenever the gap defined in~\eqref{def:adiabatic} is minimal on a trajectory. 

 \smallskip
Assume $d=2$ and consider the potential 
\[
V(x)=\begin{pmatrix}
w_1(x) & w_2(x) \\
w_2(x) & -w_1(x) 
\end{pmatrix}.
\]
Denoting by $\Pi_\pm$ the eigenprojectors of $V$, one has
\[
 V=\lambda_+\Pi_++\lambda_-\Pi_-,
\;
\lambda_\pm(x)=\pm\sqrt{w_1(x)^2+w_2(x)^2}.
\]
 Eigenvalue crossings occur on the set 
 \[
 \Upsilon=\{(x,\xi)\in\R^{4}, \; w_1(x)=w_2(x)=0\}
 \]
 that is  a submanifold of $\R^{4}$ under the assumption 
 \[ 
 {\rm Rk} \, dw_{|\Upsilon}=2.
 \]
The classical trajectories $\Phi^t_\pm$ associated with the Hamiltonian $\frac{|\xi|^2}2 +\lambda_\pm(x)$ can be continuously continued through points $(x,\xi)\in \Upsilon$ such that 
\[
dw(x)\xi:= \xi_1 \nabla w_1(x) + \xi_2 \nabla w_2(x) \not= 0_{\R^2}.
\]
The gap between the eigenvalues 
\[
g(x)=2| w(x)|
\]
is minimal along a trajectory when it passes through the hypersurface 
 \[
\Sigma=\left\{ (x,\xi)\in\R^2 ,\;\;
w(x)\cdot (dw(x)\xi)=0\right\}.
\]
This set is called {\it hoping hypersurface} in the chemical literature because switches between  modes occurs on~$\Sigma$, as we shall see now.

In order to describe the transitions, one considers an extended phase space 
\[
T_\pm^*\R^2= \R^{4}\times \{+1,-1\},
\]
 and trajectories defined on   $T_\pm^*\R^2$ as branches of smooth trajectories $(\Phi^t_\pm)_{t_i\leq t\leq t^*}$ that splits into two trajectories 
\[
  (\Phi^t_\pm)_{t^*\leq t\leq t_f}
  \;\;\mbox{and}\;\; (\Phi^t_\mp)_{t^*\leq t\leq t_f},
\]
whenever $\Phi^{t^*}_\pm\in \Sigma$. 
The initial and final times $t_i$ and $t_f$ are such that on the time interval $[t_i,t_f]$  the trajectory only  reaches $\Sigma$ at a single time $t^*$.
The probability of switching from the mode $\pm$ to the mode $\mp$ is given by the {\it Landau-Zener transition rate}
\[
T(x,\xi)=\exp \left(-\frac{\pi}{h} 
\frac {|w(x)|^2 }{|dw(x) \xi |}
\right).
\]
This generates a random walk  characterized by the probability 
$\mathbb P_{z,\ell, t}(\Gamma) $ of reaching  $\Gamma\subset T_\pm^*\R^2$ at time~$t$ starting from the point $(z,\ell)\in T_\pm^*\R^2$.
With this probability law is associated a Markov process $ \mathcal L^t_{\rm LZ} $ on the set of functions defined on  $T_\pm^*\R^2$
\[
\mathcal L_{\rm LZ}^t f(z,\ell)= \int_{T_\pm^*\R^2}
f(z',\ell')d \mathbb P _{z,\ell,t} (z',\ell').
\]
By identifying the set of observables 
\[
\mathcal A= \left\{a\in\mathbb C^{2,2},  \;a=a_+\Pi_++a_-\Pi_-\right\}
\]
 to functions on $T_\pm^*\R^2$ according to 
\[
(x,\xi,\pm1)\mapsto a_\pm(x,\xi),
\]
one extends the actions of $\mathcal L_{\rm LZ}^t$ to functions of~$\mathcal A$. 
Then, it is proved in~\cite{FL1} that under reasonable assumptions, if $\theta\in \Cinfc (\mathbb R)$ and $a\in\mathcal A$,
\begin{align*}
&\Bigl\| 
\int \theta(t) \Bigl( 
U^h(-t){\rm op}_h(a)\,U^h(t) \\
&\qquad - {\rm op}_h \left( \mathcal L_{\rm LZ}^t a\right)
\Bigr) dt
\Bigr\|_{\mathcal L(L^2(\mathbb R^2))}\leq C\,h^{1/8}.
\end{align*}
The proof of this result relies on reduction to normal forms as
in~\cite{CdV} and precise analysis of the normal forms. Finding an
optimal version of the latter estimate is open.


\subsection{Geometric aspects and  application to quantum limits}\label{sec:app_geo}

We discuss here our second application on the behavior of sequences of eigenfunctions of the Laplace-Beltrami  operator of a smooth compact manifold without boundary. 

\subsubsection{Semi-classical analysis on Riemannian manifolds}

To extend the semi-classical approach to manifolds, one needs an invariance through change of variables.
\smallskip 

With   $\kappa$ a diffeomorphism, from an open set $U$ into
$V=\kappa(U)$, is associated the  local symplectomorphism 
\[
\sigma_\kappa: \; z=(x,\xi) \mapsto \left(\kappa(x) , \, ^t d\kappa(x)^{-1} \xi\right).
\]
The map $\sigma_\kappa$ is associated with the unitary transformation
$J_\kappa$ of $L^2(\R^d)$
 \begin{equation*}
J_ \kappa f= {\rm Jac}(\kappa) ^{-\frac 12} f \circ \kappa^{-1}\in \Cinfc(V),
 \quad f \in \Cinfc(U).
\end{equation*}
There exists a constant
$C>0$ and a semi-norm $N$ such that for all $a\in \Cinfc
(V\times\R^{d})$ 
one has
\begin{align}\label{invariance}
\left\| J_\kappa^* \Op(a) J_\kappa-\Op(b)\right\|_{{\mathcal L}(L^2(\R^d))}
\leq C\, h N(a),
\end{align}
where $b = a \circ \sigma_\kappa \in \Cinfc
(U\times\R^{d})$.

This result allows one to define \hpdo on a
Riemannian manifolds $M$ through local charts. However, they are only defined at leading order (up to $O(h)$ in $\mathcal L(L^2(M))$~\cite{Zwobook}.
\smallskip

Relation~\eqref{invariance} also enlightens the geometric structure of semi-classical measures that appear as measures on the cotangent space $T^*M$, the bundle above $M$ whose fibers above $x\in M$ consists in the dual set of the tangent set $T_xM$. 
\smallskip

The semi-classical approach can also be extended in the context of (noncommutative) nilpotent graded Lie groups and nilmanifods (that are quotient of such a group by one of its co-compact subgroup), using the definition of the Fourier transform via representation theory~\cite{FF2}, or in infinite dimensional frameworks~\cite{AN,AJN}.

\subsubsection{Quantum limits}

With these tools in hand, one can consider a sequence of eigenfunctions $(\varphi_k)_{k\in\N}$ of the Laplace-Beltrami operator $-\Delta_M$ on $M$ as defined in~\eqref{eq:eigen}.
The asymptotics as $E\rightarrow +\infty$ of the  counting function
\[
N(E)=\# \{ k\in\N,\;\; E_k\leq E\}.
\]
are described by the  {\it Weyl law}~\cite{Weyl,Zwobook}
\begin{equation*}
  N(E)\sim  (2\pi)^{-d} E^{\frac d2}\, {\rm Vol}(M)\, \omega_d.
\end{equation*}
Here $d$ denotes the dimension of $M$ and $\omega_d$ is the volume of
the Euclidean unit ball in $\R^d$. 
It can be  derived from~\eqref{trace}, writing  
\[
N(E)=\| \Op(a)\|^2_{{\rm HS}(L^2(M))}
\]
for  $\Op(a )=\chi(-h^2\Delta_M)$ with $h=E^{-\frac 12}$ and $\chi\in\Cinfc(T^*M)$ approaching ${\rm 1}_{[0,1]}(\xi) {\rm 1}_M(x)$~\cite{Anantharaman_book,Zwobook}. 
\smallskip

A large literature is devoted  to the analysis of  the limit as $E\rightarrow +\infty$ of
\[
\frac 1{N(E)} \sum_{E_k\leq E}  \left| \int_M \phi(x)|\varphi_k(x) |^2 dx \right|^2,\;\; \phi\in \Con^0(M),
\]
extended to functions  $a\in\Cinfc(T^*M)$ as
\[
\nu_E(a):=\frac 1{N(E)} \sum_{E_k\leq E}\left|   \left({\rm Op}_{E^{-\frac 12}} (a) \varphi_k, \varphi_k\right) \right|^2.
\]
 The geodesic flow
 is the Hamiltonian flow associated with the symbol of the
Laplace-Beltrami operator $-\Delta_M$. Since one deduces~\eqref{eq: localization supp mu} and~\eqref{eq: transport
  mu} from equation~\eqref{eq:eigen}, properties of the geodesic flow
have consequences for the limits of
$\nu_E(a)$ as $E\rightarrow +\infty$.
\smallskip  

A flow $\Phi^t$ is said to be  {\it ergodic} if  for  Lebesgue  almost all 
$(x_0,\xi_0)\in T^*M$ and for all $a\in \Con^0(T^*M)$, 
\begin{align*}
&\frac 1T \int_0^T \Phi^t_* a(x_0,\xi_0) dt
\\
&\; \Tend{T}{+\infty}\int_{S^*M} a(x,|\xi_0| \omega) d{\rm Vol}(x)d\sigma_x(\omega).
\end{align*}
Here, $d{\rm Vol}(x)={\rm Vol}(M)^{-1} dx$ is the normalized measure on $M$ and $d\sigma_x(\omega)$ the measure on the sphere $S^*_xM$ where $(x,\omega)\in S^*M$ iff
\[
 \omega\in S_x^*M:=\{\xi\in T^*M,\;\;\|\xi \|_x=1\}.
\] 
In the formula above, $\| \cdot\|_x$ denotes the vector norm in $T_x^*M$.
The result is the following (see~\cite{Sni,CdV:85,Zel87}).

\begin{thm} \label{t:QE}
If the geodesic flow of $M$ is ergodic, then
\begin{align*}
&\lim_{E\rightarrow +\infty}\frac{1}{N(E)}
\sum_{E_k\leq E} \Biggl|\left( {\rm Op}_{E^{-\frac 12}} (a)\varphi_k,  \varphi_k\right)-\\
&\int_{S^*M} a\left(x,\left(\frac {E_k}E\right)^{1/2} \omega \right)d{\rm Vol}(x)d\sigma_x(\omega) \Biggr|^2 =0.
\end{align*}
\end{thm}

The result has an alternative equivalent version that reduces to considering the average of 
\begin{align*}
\nonumber
\mathbb L(a,E_k):=&
\Biggl|\left( {\rm Op}_{E_k^{-\frac 12}} (a)\varphi_k,  \varphi_k\right) \\
&\qquad-
\int_{S^*M} a\left(x,\omega \right)d{\rm Vol}(x)d\sigma_x(\omega) \Biggr|^2
\end{align*}
for eigenvalues $E_k$ such that $\frac E2\leq E_k\leq \frac{3E}2$. 
One then has 
\begin{align}
&
\lim_{E\rightarrow +\infty}\frac{1}{N(\frac{3E}2)-N(\frac E2)}
\sum_{\frac E2\leq E_k\leq \frac{3E}2}
\mathbb L(a,E_k)
\label{QED:equiv}
 =0.
\end{align}
At  that level,  the semi-classical
aspects are easier to see. Indeed,~\eqref{eq:eigen} shows that the family
$(\varphi_k)_{k\in\N}$ is $E_k^{-\frac 12}$-oscillating, which
motivates to adopt a semi-classical setting with $h_k = E_k^{-\frac
  12}$. The localization property~\eqref{eq: localization supp mu}
implies that the support of any semi-classical measure of
$(\varphi_k)_{k\in\N}$ is supported  in $S^*M$ and the limit
in~\eqref{QED:equiv} is the semi-classical one since $E_k\sim E$ therein. Besides, the propagation result~\eqref{eq: transport mu}   implies the invariance of the semi-classical measures of sequences $(\varphi_k)_{k\in\N}$  under the geodesic flow. When this flow is ergodic the only measure invariant under the geodesic flow is the Liouville measure, which restricts the set of quantum limits to the Liouville measure.
\smallskip 

The result is even stronger. Indeed, one of its
 consequence consists in the existence of a set $S\subset \N$ of density~$1$, meaning a set  satisfying
\[
\frac{\#\{k\in S,\, E_k\leq E\}}{N(E)} \Tend{E}{+\infty} 1
\]
such that 
\begin{align}
\nonumber
&\left( {\rm Op}_{E_k^{-\frac 12}} (a)\varphi_k,  \varphi_k.  \right) \\
\label{limit:RS}
&\;\; \Tend{k}{+\infty, \, k\in S}
\int_{S^*M} a\left(x,\omega \right)d{\rm Vol}(x)d\sigma_x(\omega) .
\end{align}
The {\it unique quantum ergodicity conjecture} of Z.~Rudnick and P.~Sarnak~\cite{Rudnick_Sarnak}  predicts that the limit in~\eqref{limit:RS} holds for the full sequence, see N.~Anantharaman's book~\cite{Anantharaman_book}.

More generally, 
 the relation between the geometry of the manifolds and the nature of quantum limits has been the subject of intensive research during the last decades (see~\cite{Hass10,ALM:16,DyatlovJin,DJN,AnNon,AM:14,MR:16, Hezari-Riviere} among others), while similar problematic arose  in other settings (for sub-Laplacian~\cite{CdVHT} and on random graphs~\cite{Nalini_Lemasson} for example).  
\smallskip


\subsection{Semi-classical methods in control theory}\label{sec:app_control}
With the notations of the preceding paragraph, 
for $n  \in \N^*$ set $J_n= \{k;\ 2^{n-2}  < \sqrt {E_k}< 2^{n+2}  \}$, 
and denote by $\mathcal F_n$ the set of functions of the form 
\begin{align*}
  u (t,x) = \sum_{k \in J_n} e^{i t \sqrt{E_k+1}} u^k \varphi_k(x), 
\end{align*}
for the coefficients $u^k \in \C$. They are solutions to the
wave equation~\eqref{eq: free wave equation} with a time-frequency $\tau \sim 2^n$. Set the
semi-classical parameter to be $h_n = 2^{-n}$.

Suppose $T>0$ and $\omega$ is an open subset of $M$.  Suppose there
exist $C>0$ and $n_0\in \N^*$ such that for all $n \geq n_0$ and all $u \in \mathcal
F_n$ one has 
\begin{align}
  \label{eq: semi-classical observation}
  \mathcal E (u) \leq C \int_0^T \Norm{1_\omega \d_t u }{L^2(M)}^2 dt, 
\end{align}
that is, observability for those frequency-localized solutions;
see~\eqref{eq: obs wave equation} where observability is defined. Then,
the wave equation \eqref{eq: free wave equation} is observable from
$\omega$ in time $T>0$, and exact controllability follows. One calls
estimate \eqref{eq: semi-classical observation} a {\it semi-classical
observability inequality}. We discuss conditions for its validity in the next Section~\ref{sec:GCC}.

The proof of this extension of observability to all solutions  
to the wave equation in \cite{Lebeau:92,BDLR:23a} makes use of the 
following unique continuation property  
\begin{align}
  \label{eq: unique continuation eigenfunction}
  - \Delta_M  \varphi = \mu \varphi \ \ \text{and} \ \ \varphi_{|\omega} =0 
 \ \ \Rightarrow \ \ \varphi =0.
\end{align}
A now classical tool to prove such a result is a Carleman estimate
that can be viewed as a sub-elliptic semi-classical estimate; these
estimates are presented in 
Section~\ref{sec:UCP}. 

In fact, the semi-classical observability estimate \eqref{eq:
  semi-classical observation} takes care of the high-frequency
component of the solutions, while the unique continuation properties
handles the remaining low frequencies.

\subsubsection{Geometric control condition}\label{sec:GCC}

To analyse observability issues, it is classical to adopt a space-time point of
  view and to work in the variables
  \[ (t,x,\tau,\xi)\in T^*(\R_t\times M_x).\] One can have intuitions
  based on geometrical optics and propagation of energy along rays to
  support this point of view. This intuition turns out to be correct
  as explained below.

The symbol of the wave operator is
\[
p(x,\tau,\xi) = -\tau^2 + g_x(\xi,\xi)
\]
in local coordinates. 
The Hamiltonian curves (rays) of the space-time Hamiltonian $p$  are
called {\it bicharacteristic curves}. Their projections on $M$ are the {\it geodesics}.

The semi-classical observability estimate~\eqref{eq: semi-classical
  observation} is proven to hold under the following property: any
\bichar reaches a point above $]0,T[ \times \omega$. This condition is
called the {\em geometric control condition} (GCC). Equivalently it
reads: any geodesic travelled at speed one
enters the observation region $\omega$ in a time less than $T$.  Then,
assuming (GCC), the proof of~\eqref{eq:
  semi-classical observation} can be carried out by
contradiction. Suppose $U = (u_n)_{n \in \N^*}$ is a sequence with
$u_n \in \mathcal F_n$ such that
\begin{align}
  \label{eq: contradiction argument}
  h_n^2 \mathcal E (u_n) \sim 1 \ \ \text{and} \ \ \int_0^T \Norm{1_\omega
  h_n\d_t u_n }{L^2(M)}^2 dt  \to 0,
\end{align}
as $n \to +\infty$. Then 
$\Norm{u_n (t,.)}{L^2(M)} \sim 1$
and associated with a subsequence
of $U$ is a semi-classical measure~$\mu$. Note that $U$ is bounded
in $L^2_{\mathrm{loc}} \big(\R_t; L^2(M)\big)$ here. The measure is
thus understood acting on functions compactly supported in the variable
$t$.
On the one hand, one has 
\begin{align*}
  h_n \Norm{\nabla_x u_n (t,.)}{L^2(M)} 
  \  \sim \  h_n \Norm{\d_t  u_n
  (t,.)}{L^2(M)} \ \sim 1,
\end{align*}
implying that $U$ is $h_n$-oscillating. One deduces that $\mu$ has positive mass. 
On the other hand, the second part of~\eqref{eq: contradiction
  argument} gives $\mu=0$ above $]0,T[\times \omega$. 

With Section~\ref{sec: semi-classical measures and PDEs}, one finds
that $\supp(\mu) \subset \Char(p)$ and $\Hp \mu=0$ in the sense of
distributions, meaning that $\mu$ is invariant along the
bicharacteristic flow. By the GCC, all \bichars enter the
region above $]0,T[\times \omega$ where $\mu$ vanishes implying that
$\mu=0$. A contradiction.
\smallskip

 Arguments are more involved in the case of a
 manifold with boundary and a wave equation formulated with a boundary
 condition, say the homogeneous Dirichlet condition. Away from the
 boundary, the measure equation $\Hp \mu=0$ holds. At the boundary, one
 can derive a measure equation that includes a source term associated
 with the semi-classical measure of the Neumann
 trace. This source term generates transport of the measure $\mu$
 along the Melrose-Sj\"ostrand generalized \bichars \cite{MS:78}. Those obey the
 laws of geometrical optics: reflection if the boundary is hit
 transversally, possible glancing and gliding if the boundary is hit
 tangentially. The GCC remains unchanged apart from exchanging
 \bichars with generalized \bichars and the observability/exact
 controlability result holds under this condition.

 The proof of wave observability with the sharp GCC condition was
 first given by C.~Bardos, G.~Lebeau, and J.~Rauch~\cite{BLR:92}, in
 the case of smooth coefficients with microlocal techniques based on
 the propagation of singularities. The use of measures was initiated
 by N.~Burq and P.~G\'erard to further explain the necessary and
 sufficient aspects of the GCC~\cite{BG:1997}. One interest of the use
 of semi-classical measures is the possibility of lowering the
 regularity of the coefficients. In \cite{BDLR:23b,BDLR:23c}, this
 regularity in pushed down to a $\Con^1$-metric on a $\Con^2$-manifold
 with boundary. Then, the Hamiltonian vector field $\Hp$ is only
 continuous. Generalized bicharacteristics exist but uniqueness is
 lost. Yet, the GCC makes sense and despite the absence of flow one
 proves that the support of the semi-classical measure $\mu$ is a
 union of generalized \bichars, allowing one to conclude the
 contradiction argument as above.

\subsubsection{Carleman estimates as sub-elliptic semi-classical
  estimates and unique continuation}\label{sec:UCP}

For a second-order elliptic operator $P$, a
{\it Carleman estimate} takes the
form, for some $C>0$, 
\begin{multline}
  \label{eq: Carleman}
  h^{1/2} \big( \Norm{e^{\varphi/h} u}{L^2}
  + \Norm{e^{\varphi/h} h \nabla_x u}{L^2}\big) \\
  \leq C \Norm{h^2 e^{\varphi/h} P u}{L^2}, 
\end{multline}
for $u$ smooth with compact support. 
The inequality holds if the function $\varphi$,  called the {\it weight
function}, is well chosen and if  $0< h \leq h_0$,  for $h_0$
\suff small. 
For $x_0 \in \R^d$,  a possible choice of weight function is 
\begin{align*}
  \varphi(x) = \exp(-\gamma |x-x_0|). 
\end{align*}
Given any $c_0>0$, estimate \eqref{eq: Carleman} holds
for functions $u$ supported in the annulus $0< c_0 \leq |x-x_0| \leq
 4 c_0$ if  $\gamma>0$ is chosen large
\cite{Hoermander:V1,LRL:12,LRLR:V1}. For $a>0$, set $B_a  = \{ x\in\R^d,\, |x-x_0| \leq
 a \}$
 and $B' = B_{4c_0}$.
  Then, for $c_0< r <  2 c_0$ one
deduces from~\eqref{eq: Carleman} the existence of $C>0$ and $\delta \in ]0,1[$ such that 
\begin{align*}
  \Norm{u}{H^1(B_{2r})} 
  \leq C \Norm{u}{H^1(B')}^{1-\delta} 
  \big(
  \Norm{Pu}{L^2(B')} + \Norm{u}{H^1(B_r)} \big)^\delta.
\end{align*}
This is a quantification of the unique continuation property: if $u=0$
in $B_r$ and $Pu=0$ then $u=0$ in $B_{2r}$. Applied to $P = -\Delta_M
-\mu$, one obtains \eqref{eq: unique continuation eigenfunction}.

Estimate \eqref{eq: Carleman} is equivalent to 
\begin{align}
  \label{eq: Carleman-2}
  h^{1/2} \big( \Norm{v}{L^2}
  + \Norm{h \nabla_x v}{L^2}\big) 
  \lesssim \Norm{P_\varphi v}{L^2}, 
\end{align}
where $P_\varphi = h^2 e^{\varphi/h} P e^{-\varphi/h}$ is a
semi-classical differential operator. This operator
fails to be elliptic in general, yet the weight function $\varphi$ is
chosen so that the  following property holds for  $p_\varphi$, the symbol of $P_\varphi$, 
\begin{align*}
  p_\varphi(x,\xi) =0 \ \ \Rightarrow \ \ 
  \{ \Re p_\varphi, \Im p_\varphi\} = \frac{1}{2i} \{ \ovl{p_\varphi},
  p_\varphi\} >0,
\end{align*}
that is, a subellipticity property. This explains the
factor $h^{1/2}$ on the left-hand side of  \eqref{eq: Carleman} that
expresses a half-derivative loss as compared to an elliptic
estimate. One proves
\begin{align*}
  \nu |p_\varphi|^2 (x,\xi) 
  + \{ \Re p_\varphi, \Im p_\varphi\}  (x,\xi) 
  \geq C(1 + |\xi|^4),
\end{align*}
 and \eqref{eq: Carleman-2} follows from the sharp G{\aa}rding inequality
 \eqref{eq: gaarding}.

Carleman inqualities can be derived for other types of operators and
under fine geometrical properties between the operator $P$ and the
weight function $\varphi$, so-called pseudo-convexity conditions. The
reader is referred to \cite{Lerner:19} for an exposition.

\section{Concluding remark}

In the 1970s, people from the microlocal community started to show a
strong interest in semi-classical analysis. Let us mention some of the
first contributions on the domain by V. Maslov~\cite{maslov} and
A.~Voros~\cite{voros}.
 Fifty years later, the
theory has grown considerably.  Taking the correspondence principle as a guideline, the authors aimed to
show how vast the field of applications of semi-classical is today. The authors hope they have managed to pass
on their interest and enthousiasm for semi-classical analysis through
striking results. The
different examples presented here reflect the mathematical tastes of the
authors and should not be thought as exhaustive.
They
recommend the reading of the section devoted to semi-classical
analysis by B.~Helffer in the previous edition of this
encyclopedia, his commented bibliography~\cite{helf_bib};  the books by M.~Zworski~\cite{Zwobook} and by M.~Dimassi and J.~Sj\"ostrandt~\cite{DimassiSjostrand} will be useful  for
both junior and confirmed researchers.

\section*{Acknowledgements}
The authors warmly thank Maher Zerzeri for his remarks and comments, Caroline Lasser and Fabricio Macia for inspiring discussions.


\def\cprime{$'$}

\noindent {\bf Adresses}\\ {CFK - Larema, UMR  6093,
universit\'e d'Angers, CNRS,
clotilde.fermanian@univ-angers.fr\\
JLR - LAGA, UMR 7339,
université Sorbonne Paris Nord, CNRS, 
jerome.lerousseau@univ-paris13.fr.}

\end{document}

%% file: macros.tex
\let \Re \relax
\DeclareMathOperator{\Re}{Re}
\let \Im \relax
\DeclareMathOperator{\Im}{Im}

\newcommand{\Op}{{\rm Op}_h}

\def\Tend#1#2{\mathop{\longrightarrow}\limits_{#1\rightarrow#2}}

\newcommand{\mb}[1]{\ensuremath{\mathbb{#1}}}
\newcommand{\N}{{\mb{N}}}

\newcommand{\R}{{\mb{R}}}
\newcommand{\C}{{\mb{C}}}

\renewcommand{\S}{\ensuremath{\mathscr S}}

\newcommand{\tF}{\ensuremath{\tilde{F}}}

\newcommand{\transp}[1]{\ensuremath{\leftidx{^t}{\!#1}{}}}

\DeclareMathOperator{\Hamiltonian}{H}
\newcommand{\Hp}{\Hamiltonian_p}

\newcommand{\Con}{\ensuremath{\mathscr C}}

\newcommand{\Cinfc}{\ensuremath{\mathscr C^\infty_{c}}}

\renewcommand{\d}{\ensuremath{\partial}}

\newcommand{\inp}[2]{(#1, #2)} 
\newcommand{\biginp}[2]{\big(#1, #2 \big)}

\newcommand{\dup}[2]{\langle #1, #2 \rangle}

\newcommand{\ovl}[1]{\overline{#1}}

\newcommand{\jpsi}[1]{\langle #1 \rangle_{\mbox{\tiny\!$\psi$}}}
\newcommand{\bjpsi}[1]{\big\langle #1 \big\rangle_{\mbox{\tiny\!\!$\psi$}}}
\newcommand{\Norm}[2]{{\| #1 \|}_{#2}}

\DeclareMathOperator{\Char}{Char}

\DeclareMathOperator{\supp}{supp}

\newcommand{\bichar}{bicharacteristic\xspace}
\newcommand{\bichars}{bicharacteristics\xspace}

\newcommand{\eg}{e.g.\@\xspace}

\newcommand{\hpdo}{$h$-$\psi$do\xspace}

\newcommand{\suff}{sufficiently\xspace}

\newcommand{\bld}[1]{\mbox{\boldmath $#1$}}